\documentclass[12pt,a4paper]{amsart}
\usepackage{amssymb,amsmath}
\usepackage{times}

\numberwithin{equation}{section}
\newtheorem{theorem}{Theorem}[section]
\newtheorem{lemma}[theorem]{Lemma}

\newtheorem{corollary}[theorem]{Corollary}

\theoremstyle{definition}

\newtheorem{remark}[theorem]{Remark}

\newtheorem{example}[theorem]{Example}

\newtheorem{acknowledgement}{Acknowledgement}

\newcommand\Proj{\operatorname{Proj}}

\newcommand\height{\operatorname{height}}
\newcommand\depth{\operatorname{depth}}
\newcommand\codim{\operatorname{codim}}
\newcommand\codepth{\operatorname{codepth}}
\newcommand\reg{\operatorname{reg}}
\newcommand\rank{\operatorname{rank}}
\newcommand\Tor{\operatorname{Tor}}

\newcommand\Sec{\operatorname{Sec}}

\begin{document}

\author[M.~Brodmann, P.~Schenzel]{Markus Brodmann,  Peter Schenzel}
\title[Varieties of almost minimal degree]
{On varieties of almost minimal degree in small codimension}

\thanks{The second named author was partially supported by Swiss National Science
           Foundation (Project No. 20 - 103491/1 )}

\address{Universit\"at Z\"urich, Institut f\"ur Mathematik,
Winterthurer Str. 190, CH -- 8057 Z\"urich,
Switzerland} \email{brodmann@math.unizh.ch}
\address{Martin-Luther-Universit\"at Halle-Wittenberg,
Fachbereich Mathematik und Informatik, D -- 06 099 Halle (Saale),
Germany} \email{schenzel@mathematik.uni-halle.de}

\begin{abstract} 
The present research grew out of the authors' joint work
\cite{BS}. It continues the study of the structure of
projective varieties of almost minimal degree, focussing to the 
case of small codimension. In particular, we give a complete list of 
all occuring Betti diagrams in the cases where $\codim X \leq 4.$ 
\end{abstract}

\maketitle

\section{Introduction}
Let $X \subset \mathbb P^r_K$ denote an irreducible and reduced
projective variety over an algebraically closed field $K.$ We 
always assume that $X$ is non-degenerate, that is not
contained in a hyperplane. Then, the degree and the codimension of
$X$ satisfy the inequality $\deg X \geq \codim X +1$ (cf. for
instance \cite{H}). Varieties for which equality holds are called
called {\sl varieties of minimal degree}. These varieties are completely
classified (cf. for instance \cite{EH} and \cite{H}). In
particular they are arithmetically Cohen-Macaulay and have a
linear minimal free resolution. In particular, their Betti numbers are
explicitly known.

In  case $\deg X = \codim X +2,$ the variety $X$ is called a
{\sl variety of almost minimal degree}. Here one has a much
greater variety of possible Betti numbers. The investigation of 
homological properties of varieties of almost
minimal degree was initiated by Hoa, St\"uckrad, and Vogel (cf.
\cite{HSV}). We refer also to \cite{BS} and \cite{Na} for certain
improvements of their results. In
particular the Castelnuovo-Mumford regularity of a variety $X$ of
almost minimal degree satisfies $\reg X \leq 2$ (cf. \cite{E} for
the definition of the Castelnuovo-Mumford regularity).

In the framework of
polarized varieties of $\Delta$-genus 1, Fujita (cf. \cite{Fu1}
and \cite{Fu2}) provides a satisfactory description of varieties
of almost minimal degree. The study of varieties of almost minimal
degree is pursued by the authors (cf. \cite{BS}) from the
arithmetic point of view. It turns out that $X
\subset \mathbb P^r_K$ is either an arithmetically normal Del
Pezzo variety or a proper projection of a variety of minimal
degree. By a proper projection of a variety $Z \subset \mathbb P^{r+1}_K$ 
we always mean a projection from a point $P \in \mathbb P^{r+1}_K \setminus 
Z.$ See also \ref{L 3.1} for the precise statement.

The aim of the present paper is to investigate varieties of
almost minimal degree and of low codimension, in particular their
Betti diagrams. More precisely, we describe the structure of the
minimal free resolution of a variety $X$ of almost minimal degree
of $\codim X \leq 4$ by listing all possible Betti diagrams. Let
us recall that the structure of arithmetically Cohen-Macaulay
resp. Gorenstein varieties in codimension 2 resp. 3 is known by
the Theorems of Hilbert-Burch resp. Buchsbaum-Eisenbud (cf.
\cite{E}). So, we need not to discuss these cases in detail any
more. As the Betti diagram, the degree and the codimension are
not affected if $X$ is replaced by a cone over $X$, we shall
assume that $X$ is not a cone. The most surprising fact is, that
the dimension of $X$ is always $\leq 6$ (cf. Section 2 for the 
precise statements). Our main
technical tool is a result shown by the authors in \cite{BS},
which says that apart from an exceptional case, (that is the
generic projection of the Veronese surface in $\mathbb P^5_K$) any
non-arithmetically normal (and in particular non-arithmetically
Cohen-Macaulay) variety of almost minimal degree $X \subset
\mathbb P^r_K$ (which is not a cone) is contained in a variety of
minimal degree $Y \subset \mathbb P^r_K$ such that $\codim(X,Y) =
1.$

\begin{acknowledgement} The authors thank Uwe Nagel for making his preprint
\cite{Na} available. They also thank for the referee's very valuable 
suggestions.
\end{acknowledgement}

\section{Main Results}
Let $X \subset \mathbb P^r_K$ denote a non-degenerate variety of
almost minimal degree, hence an integral closed subscheme with
$\deg X = \codim X +2$ not contained in a hyperplane $\mathbb
P^{r-1}_K \subset \mathbb P^r_K.$ We use the abbreviations $\dim X
= d$ and $ \codim X = c.$ Let $S= K[x_0,\ldots,x_r]$ denote the
polynomial ring in $r+1$ variables, so that  $\mathbb P^r_K =
\Proj(S).$ Let $A_X = S/I_X$ denote the homogeneous coordinate
ring of $X,$ where $I_X \subset S$ is the defining ideal of $X.$
The $\codepth$ of $A_X$ is defined as the difference $\codepth
A_X := \dim A_X - \depth A_X,$ where $\depth A_X$ denotes the
depth of $A_X.$

For the notion of Betti diagrams we follow the suggestion of
Eisenbud (cf. \cite{E}). That is, in a diagram the number in the
$i$-th column and the $j$-th row is
\[
\dim_K \Tor_i^S(K,S/I_X)_{i+j}, \; i \geq 1.
\]
Outside the range of the diagram all the entries are understood to
be zero. Our main results are.

\begin{theorem} \label{T 2.1} $(\codim X = 2)$ Let $X \subset
\mathbb P^r_K$ be a non-degenerate variety of degree 4 and
codimension 2 which is not a cone. Then $X$ is of one of the
following types:
\begin{itemize}
\item[(a)] $X$ is a complete intersection cut out by two quadrics.
\item[(b)] $\dim X \leq 4$ and the Betti diagram of $I_X$ has the
form
\begin{center}
\begin{tabular}{|c|ccc|}
  \hline
   & 1 & 2 & 3 \\
  \hline
  1 & 1 & 0 & 0 \\
  2 & 3 & 4 & 1 \\
  \hline
\end{tabular}\; .
\end{center}
\item[(c)] (The exceptional case) $X$ is a generic projection of
the Veronese surface $F \subset \mathbb P^5_K$ and the Betti
diagram of $I_X$ has the shape
\begin{center}
\begin{tabular}{|c|cccc|}
  \hline
    & 1 & 2 & 3 & 4 \\
  \hline
  1 & 0 &  0 & 0 & 0 \\
  2 & 7 & 10 & 5 & 1 \\
  \hline
\end{tabular}\; .
\end{center}
\end{itemize}
Moreover for any $1 \leq d \leq 4$ there are examples as mentioned
in (b) such that $\dim X = d.$
\end{theorem}

In the case where $X$ is a Cohen-Macaulay variety, Theorem \ref{T
2.1} (b) has been shown by Nagel (cf. \cite{Na}). Under this
additional assumption one has $\dim X \leq 2.$ Theorem \ref{T 2.1}
grew out of our aim to understand Nagel's arguments. Our approach
enables us to investigate the cases of codimension three and four
as well.

\begin{theorem} \label{T 2.2} $(\codim X = 3)$ Let $X \subset
\mathbb P^r_K$ denote a non-degenerate variety of degree 5 and
codimension 3 which is not a cone. Then the
following cases may occur:
\begin{itemize}
\item[(a)] $\dim X \leq 6$ and $X$ is the Pfaffian variety defined by the five
Pfaffians of a skew symmetric $5 \times 5$ matrix of linear forms. Its minimal
free resolution is given by the Buchsbaum-Eisenbud complex.
\item[(b)] $\dim X \leq 4$ and $\codepth A_X = 1.$ The Betti
diagram of the defining ideal $I_X$ has the following form
\begin{center}
\begin{tabular}{|c|cccc|}
  \hline
    & 1 & 2 & 3 & 4 \\
  \hline
  1 & 4 & 2 & 0 & 0 \\
  2 & 1 & 6 & 5 & 1 \\
  \hline
\end{tabular}\; .
\end{center}
\item[(c)] $\dim X \leq 5$ and $\codepth A_X = 2.$ The Betti
diagram of $I_X$ has the following shape
\begin{center}
\begin{tabular}{|c|ccccc|}
  \hline
    & 1 & 2 & 3 & 4 & 5 \\
  \hline
  1 & 3 & 2 & 0 & 0 & 0 \\
  2 & 6 & 16 & 15 & 6 & 1  \\
  \hline
\end{tabular}\; .
\end{center}
\end{itemize}
For any of the dimensions and codepths admitted in (a), (b) and (c)
resp. there are examples of varieties of almost minimal degree.
\end{theorem}

The next result concerns the case where $X$ is of codimension 4. A new
phenomenon occurs in this situation. Namely, for the same codepth two
different Betti diagrams may occur.

\begin{theorem} \label{T 2.4} $(\codim X = 4)$ Let $X \subset
\mathbb P^r_K$ denote a non-degenerate variety of degree 6 and
codimension 4 which is not a cone. Then the following four cases 
may occur:
\begin{itemize}
\item[(a)] $\dim X \leq 4$ and $X$ is arithmetically Gorenstein. Its minimal free
resolution has the following form
\[
0 \to S(-6) \to S^9(-4) \to S^{16}(-3) \to S^9(-2) \to I_X \to 0.
\]

\item[(b)] $\dim X \leq 4$ and $\codepth A_X = 1.$ The Betti
diagram of the defining ideal $I_X$ has one of the following two
forms:
\begin{center}
\begin{tabular}{|c|ccccc|}
  \hline
    & 1 & 2 & 3 & 4 & 5 \\
  \hline
  1 & 8 & 12 & 3 & 0 & 0 \\
  2 & 1 & 4 & 10 & 6 & 1 \\
  \hline
\end{tabular}
\mbox{ resp. }
\begin{tabular}{|c|ccccc|}
  \hline
    & 1 & 2 & 3 & 4 & 5 \\
  \hline
  1 & 8 & 11 & 3 & 0 & 0 \\
  2 & 0 & 4 & 10 & 6 & 1 \\
  \hline
\end{tabular}\; .
\end{center}
\item[(c)] $\dim X \leq 5$ and $\codepth A_X = 2.$ The Betti
diagram of $I_X$ is of the form:
\begin{center}
\begin{tabular}{|c|cccccc|}
  \hline
    & 1 & 2 & 3 & 4 & 5 & 6\\
  \hline
  1 & 7 & 8 & 3 & 0 & 0 & 0\\
  2 & 3 & 19 & 30 & 21 & 7 & 1 \\
  \hline
\end{tabular}\; .
\end{center}
\item[(d)] $\dim X \leq 6$ and $\codepth A_X = 3.$ The Betti
diagram of $I_X$ has the following shape:
\begin{center}
\begin{tabular}{|c|ccccccc|}
  \hline
    & 1 & 2 & 3 & 4 & 5 & 6 & 7\\
  \hline
  1 & 6 & 8 & 3 & 0 & 0 & 0 & 0\\
  2 & 10 & 40 & 65 & 56 & 28 & 8 &1\\
  \hline
\end{tabular}\; .
\end{center}
\end{itemize}
For any of the dimensions, codepths and Betti diagrams 
admitted in (a), (b), (c) and
(d) resp. there are examples of varieties of almost minimal
degree.
\end{theorem}

The final result of this note concerns varieties of almost minimal
degree and of codimension less than half the embedding dimension.
First of all note (cf. Example \ref{E 4.1}) that in codimension 6
there are varieties of almost minimal degree having the same 
$\codepth$ but with rather different
Betti diagrams. That is, the corresponding statements to Theorem
\ref{T 2.1}, \ref{T 2.2} and \ref{T 2.4} are not true in hihgher
codimension.

\begin{corollary} \label{C 2.3} Let $X \subset \mathbb P^r_K$ be a
non-degenerate variety of almost minimal degree which is not a
cone. Suppose that $\dim X > \codim X
+ 2$ and $\codim X \geq 3.$ Then $\codim X = 3$ and $X$ is 
arithmetically Gorenstein.
Therefore $\dim X \leq 6$ and $X$ is defined by the five Pfaffians of a skew symmetric
$5 \times 5$ matrix of linear forms.
\end{corollary}

Moreover, in the case where $X$ is a variety of almost minimal
degree and $\dim X \leq \codim X +2,$ that is when $X$ is not
necessary a Del Pezzo variety there are estimates for the Betti
numbers and their vanishing (cf. \ref{L 3.2} for the details).

\section{Outline of the Proofs}
Let $X \subset \mathbb P^r_K$ denote a non-degenerate reduced
irreducible variety of almost minimal degree. By the work of the
authors (cf. \cite{BS}) it follows that $X$ is either a normal Del
Pezzo variety -- in this case $X$ is arithmetically Gorenstein --
or a projection of a variety of minimal degree.

First, we consider the case of non-arithmetically normal 
varieties of almost minimal degree. In this situation we have the
following characterization in which $\Sec_P(Z)$ denotes the secant
cone of a projective variety $Z$ with respect to a point $P$ in
the ambient space.

\begin{lemma} \label{L 3.1} Let $X \subset \mathbb P^r_K$ denote a
non-degenerate reduced irreducible variety which is not a cone. Let $1 \leq t \leq
\dim X + 1 =: d + 1.$ Then the following conditions are
equivalent:
\begin{itemize}
\item[(i)] $X$ is not arithmetically normal, $\deg X = \codim X + 2$ and $\depth A_X = t.$
\item[(ii)] $X$ is the projection of a variety $Z \subset \mathbb
P^{r+1}_K$ of minimal degree from a point $P \in \mathbb P^{r+1}_K
\setminus Z$ such that $\dim \Sec_P(Z) = t-1.$
\end{itemize}
Moreover, $1 \leq \depth A_X \leq 4.$
\end{lemma}

\begin{proof} Cf. \cite[Theorem 1.1 and Corollary 7.6]{BS}. \end{proof}

In view of Lemma \ref{L 3.1} there is some need for information
about varieties of minimal degree in order to understand varieties
of almost minimal degree. A variety of minimal degree $Z \subset
\mathbb P^s_K$ is either
\begin{itemize}
\item a quadric hypersurface,

\item a (cone over a) Veronese surface in $\mathbb P^5_K,$ or

\item a (cone over a) smooth rational normal scroll
\end{itemize}
(cf. \cite{EH} and \cite{H} for the details and the history of
this classification).

Next we recall a few basic facts about rational normal scrolls
(cf. also \cite{H}). Let 
\[
T = K[x_{10},\ldots,x_{1a_1},x_{20},\ldots,x_{2a_2},\ldots,x_{k0},
\ldots,x_{ka_k}]
\]
be the polynomial ring. A (cone over a) rational normal scroll
$S(a_1,\ldots,a_k)$ is defined as the "rank two subscheme in
$\mathbb P^s_K = \Proj(T)$" defined by the matrix
\[
M = 
\left(%
\begin{array}{ccccccccccccc}
  x_{10} & \ldots & x_{1a_1-1} & \vdots & x_{20} & \ldots & x_{2a_2-1} & \vdots & \ldots & \vdots & x_{k0} & \ldots & x_{ka_k-1} \\
  x_{11} & \ldots & x_{1a_1}   & \vdots & x_{21} & \ldots & x_{2a_2}   & \vdots & \ldots & \vdots & x_{k1} & \ldots & x_{ka_k} \\
\end{array}%
\right)
\]
with $s = k-1 + \sum_{i=1}^k a_i.$ It is well known that $\dim
S(a_1,\ldots,a_k) = k$ and therefore
\[
\deg S(a_1,\ldots,a_k) = \codim S(a_1,\ldots,a_k) +1 =
\sum_{i=1}^k a_i,
\]
so that $S(a_1,\ldots,a_k)$ is a variety of minimal degree. Keep
in mind that $S(a_1,\ldots,a_k)$ is a proper cone if and only if
$a_i = 0$ for some $i \in \{1,\ldots,k\},$ that is, if and only if 
there are indeterminates which do not occur in $M.$

Moreover we need some information about the Hilbert series
$F(\lambda, A_X)$ of the graded $K$-algebra $A_X.$ Let us recall
that the Hilbert series of a graded $K$-algebra $A$ is the formal
power series defined by
\[
F(\lambda, A) = \sum_{i \geq 0} (\dim_K A_i) \lambda^i.
\]
The Hilbert series of a variety of almost minimal degree may be
described as follows.

\begin{lemma} \label{L 3.4} Let $X \subset \mathbb P^r_K$ denote a
variety of almost minimal degree. Put $q = \codepth A_X.$ Then
\[
F(\lambda, A_X) = \frac{1}{(1-\lambda)^{d+1}} (1 + (c+1)\lambda -
\lambda (1-\lambda)^{q+1}),
\]
where $c = \codim X$ and $d = \dim X.$ Furthermore $\dim_K (I_X)_2
= \binom{c+1}{2} -q -1.$
\end{lemma}

\begin{proof} Cf. \cite[Corollary 4.4]{BS}. \end{proof}

As a consequence of Lemma \ref{L 3.1} the authors (cf. \cite{BS})
derived some information about the Betti numbers of $I_X$ for
certain varieties $X \subset \mathbb P^r_K$ of almost minimal
degree.

\begin{lemma} \label{L 3.2} Let $X \subset \mathbb P^r_K$ be a
variety of almost minimal degree which is not arithmetically
Cohen-Macaulay. Suppose that $X$ is not a generic projection of (a
cone over) the Veronese surface $F \subset \mathbb P^5_K.$  Then
there exists a variety of minimal degree $Y \subset \mathbb P^r_K$
such that $X \subset Y$ and $\codim(X,Y) = 1.$ Moreover
\[
\Tor_i^S(k, A_X) \simeq k^{u_i}(-i-1) \oplus k^{v_i}(-i-2) \text{
for } 1 \leq i \leq c+q,
\]
where $c = \codim X, q = \codepth A_X$ and
\begin{itemize}
\item[(a)]
\begin{itemize}
\item[] $u_1 = \binom{c+1}{2} -q-1,$

\item[] $i \binom{c}{i+1} \leq u_i \leq (c+1)\binom{c}{i} -
\binom{c}{i+1}, \text{if} \quad 1 < i < c-q,$

\item[] $u_i = i \binom{c}{i+1}, \text{if} \quad c-q \leq i < c,$

\item[] $u_i = 0, \text{if} \quad c \leq i \leq c+q.$

\end{itemize}
\item[(b)]
\begin{itemize}
\item[] $\max \{ 0, \binom{c+q-1}{i+1} - (i+2)\binom{c}{i+1}\}
\leq v_i \leq \binom{c+q+1}{i+1}, \text{if} \quad 1 \leq i <
c-q-1,$

\item[] $v_i = \binom{c+q+1}{i+1} - (i+2)\binom{c}{i+1}, \text{if}
\quad \max\{ 1,c-q-1\} \leq i < c,$

\item[] $v_i = \binom{c+q+1}{i+1}, \text{if} \quad c \leq i \leq
c+q.$
\end{itemize}
\end{itemize}
In addition $v_i - u_{i+1} = \binom{c+q+1}{i+1}
-(c+1)\binom{c}{i+1} + \binom{c}{i+2}$ for all $1 \leq i < c.$
\end{lemma}

\begin{proof} Cf. \cite[Theorem 1.1 and Theorem 8.3]{BS}. \end{proof}

In the particular case where $\dim X = 1$ the statement about the
Betti numbers has been shown independently by Nagel (cf.
\cite{Na}).

In the following remark, we  add a comment concerning the
"exceptional cases" of a generic projection of (a cone over) the
Veronese surface $F \subset \mathbb P^5_K$ and of an arithmetically 
Cohen-Macaulay variety.

\begin{remark} \label{R 3.3} A) (The exceptional case) Let $F
\subset \mathbb P^5_K$ be the Veronese surface defined by the
$2\times 2$-minors of the symmetric matrix
\[
M =
\left(%
\begin{array}{ccc}
  x_0 & x_1 & x_2 \\
  x_1 & x_3 & x_4 \\
  x_2 & x_4 & x_5 \\
\end{array}%
\right).
\]
Let $P \in \mathbb P^5_K \setminus F$ denote a point. Suppose that
$\rank M \hspace{-.2cm}\mid_P = 3,$ that is the case of a generic
point. Remember that $\det M = 0$ defines the secant variety of
$F.$ Then the projection of $F$ from $P$ defines a surface $X
\subset \mathbb P^4_K$ of almost minimal degree and $\depth A_X =
1.$ The surface $X$ is cut out by seven cubics (cf. \ref{L 3.4}),
so that it is not contained in a variety of minimal degree. \\
B) (The arithmetically normal case) Let $X \subset \mathbb
P^r_K$ denote an arithmetically normal variety of almost minimal
degree. Then $X$ is not a birational projection of a scroll and 
hence a maximal Del Pezzo variety (cf. \cite[Theorem 1.2]{BS}). In 
particular $X$ is arithmetically Gorenstein. If in addition $\codim X \geq 3,$ 
Fujita's classification of normal maximal Del Pezzo varieties yields 
$\dim X \leq 4$ (cf. \cite[(8.11), (9.17)]{Fu3}). If $\codim X \geq 4,$ 
the same classification shows that $\dim X \leq 4.$ \\
C) If $X$ is an arithmetical normal Del Pezzo variety it is 
in general not a one codimensional subvariety of a variety of minimal degree. 
Namely, 
let $X \subset \mathbb P^9_K$ be the smooth codimension three
variety cut out by the $4 \times 4$-Pfaffians of a generic
skew-symmetric $5\times 5$-matrix of linear forms. Then $X$ is not contained in a
variety $Y$ of minimal degree such that $\codim (X, Y) = 1$ (cf.
\cite{BS} for the details).
\end{remark}

\noindent {\em Proofs}. After these preparations, we now come to
the proofs of our Theorems. The proof of statement
\ref{T 2.1} (a) is easy. Now let us consider the statements \ref{T 2.2} (a) 
and \ref{T 2.4} (a). In both cases $X$ is arithmetically Gorenstein. 
If $X$ is a birational projection of a scroll, we have $\dim A_X = 
\depth A_X \leq 4$ (cf. Lemma \ref{L 3.1}). Otherwise $X$ is arithmetically 
normal (cf. \cite[Theorem 1.2]{BS}). So, by \ref{R 3.3} B) we have $\dim X \leq 6$ 
if $\codim X \geq 3$ and $\dim X \leq 4$ if $\codim X \geq 4.$ This gives 
us the dimension estimates. The statements on the minimal free resolutions 
now follow from well known results. In fact the structure of the minimal free
resolution of $I_X$ given in Theorem \ref{T 2.4} (a) is a
consequence of \cite[Theorem B]{S}.

Next we prove statement (c) of Theorem \ref{T 2.1}. To this end let $X
\subset \mathbb P^4_K$ be a generic projection of the Veronese
surface $F \subset \mathbb P^5_K.$ Then $\dim X = 2, \codim X = 2$
and $\depth A_X = 1.$ As seen above, $I_X$ does not contain any
quadric. Therefore $I_X$ has a linear resolution. Remember that
$\reg I_X = 3$ (cf. \cite{BS}). A computation with the aid of the
Hilbert series (cf. \ref{L 3.4}) gives the structure of the Betti
diagram of 2.1 (c).

For all other statements of Theorems \ref{T 2.1}, \ref{T 2.2} and
\ref{T 2.4} we may assume that the variety of almost minimal
degree $X \subset \mathbb P^r_K$ is not arithmetically
Cohen-Macaulay and not a projection of the Veronese surface $F
\subset \mathbb P^5_K.$ So, $X$ is contained in a variety of
minimal degree $Y$ of codimension $c-1$ (cf. \ref{L 3.2}). That
is, $Y$ is defined as the zero locus of $\binom{c}{2}$ quadrics.
On the other hand the defining ideal $I_X$ is generated by
$\binom{c+1}{2} -q -1$ quadrics (cf. \ref{L 3.4}). This implies
that
\[
\dim_K (I_X/I_Y)_2 = c -q -1 \geq 0,
\]
where $c = \codim X, q = \codepth A_X \geq 1.$ Considering all
possibilities that arise for $c = 2, 3, 4$ it follows that, with
the exception of the case in which $c = 4, q = 1,$ Lemma \ref{L
3.2} furnishes the corresponding Betti diagrams.

The particular case where $c = 4, q = 1,$ yields the following shape
of the Betti diagram
\begin{center}
\begin{tabular}{|c|ccccc|}
  \hline
    & 1 & 2 & 3 & 4 & 5 \\
  \hline
  1 & 8 & $u_2$ & 3 & 0 & 0 \\
  2 & $v_1$ & 4 & 10 & 6 & 1 \\
  \hline
\end{tabular}
\end{center}
with $v_1 - u_2 = -11.$ In order to finish the proof we observe
that $v_1 \leq 1$ (cf. Lemma \ref{L 3.5}).

To complete the proof of Theorems \ref{T 2.1}, \ref{T 2.2} and
\ref{T 2.4} we have to prove the stated constraints on the
occuring dimensions and codepths. To do so, we may assume that $X$
is not arithmetically Cohen-Macaulay. Therefore by Lemma \ref{L
3.1} $X$ is the projection of a variety of minimal degree $Z
\subset \mathbb P^{r+1}_K$ from a point $P \in \mathbb P^{r+1}_K
\setminus Z$ such that $\dim \Sec_P(Z) = t-1,$ where $t = \depth
A_X.$

Next let us analyze this situation in more detail. To this end let
$Z = S(a_1,\ldots,a_k)$ for certain integers $a_i, i =
1,\ldots,k.$ Then it follows that
\[
r = k-2 + \sum_{i=1}^k a_i \mbox{ and } c + 2 = \sum_{i=1}^k a_i,
\]
where $c = \codim X.$

As $X$ is not a cone, $Z$ cannot be a cone over a rational normal
scroll. Therefore $\min \{ a_i : i = 1,\ldots, k\} \geq 1.$ So,
for a given codimension $c$ we have to investigate all the
possible partitions
\[
c+2 = \sum_{i=1}^k a_i, \text{ with } k \geq 1 \text{ and } a_1
\geq a_2 \geq \ldots \geq a_k \geq 1.
\]
For $c=2$ we thus get the following possible types for the rational
normal scroll $Z:$
\[
\begin{array}{ccccccc}
  k & a_1 & a_2 & a_3 & a_4  & r+1 & \dim X \\
  \hline
  1 & 4 &   &   &   &  4 & 1 \\
  2 & 3 & 1 &   &   &  5 & 2 \\
  2 & 2 & 2 &   &   &  5 & 2 \\
  3 & 2 & 1 & 1 &   &  6 & 3 \\
  4 & 1 & 1 & 1 & 1 &  7 & 4 \\
\end{array}
\]
This proves already that $\dim X \leq 4.$

Next, we discuss the case in which the codimension equals $3.$
Here, there are the following possibilities for the type of the
scroll $Z$:
\[
\begin{array}{cccccccc}
  k & a_1 & a_2 & a_3 & a_4 & a_5 & r+1 & \dim X \\
  \hline
  1 & 5 &   &   &   &   & 5 & 1 \\
  2 & 4 & 1 &   &   &   & 6 & 2 \\
  2 & 3 & 2 &   &   &   & 6 & 2 \\
  3 & 3 & 1 & 1 &   &   & 7 & 3 \\
  3 & 2 & 2 & 1 &   &   & 7 & 3 \\
  4 & 2 & 1 & 1 & 1 &   & 8 & 4 \\
  5 & 1 & 1 & 1 & 1 & 1 & 9 & 5 \\
\end{array}
\]
Therefore $\dim X \leq 5.$ We already know that $u_1 \geq 3.$ Now,
on use of Lemma \ref{L 3.2} we easily get the requested constraints in
Theorem \ref{T 2.2}.

Finally, let $\codim X = 4.$ Then, as above there is the following
list of possible types for the scroll $Z:$
\[
\begin{array}{ccccccccc}
  k & a_1 & a_2 & a_3 & a_4 & a_5 & a_6 & r+1 & \dim X \\
  \hline
  1 & 6 &   &   &   & &  & 6 & 1 \\
  2 & 5 & 1 &   &   & &  & 7 & 2 \\
  2 & 4 & 2 &   &   & &  & 7 & 2 \\
  2 & 3 & 3 &   &   & &  & 7 & 2 \\
  3 & 4 & 1 & 1 &   & &  & 8 & 3 \\
  3 & 3 & 2 & 1 &   & &  & 8 & 3 \\
  3 & 2 & 2 & 2 &   & &  & 8 & 3 \\
  4 & 3 & 1 & 1 & 1 & &  & 9 & 4 \\
  4 & 2 & 2 & 1 & 1 & &  & 9 & 4 \\
  5 & 2 & 1 & 1 & 1 &1&  & 10 & 5 \\
  6 & 1 & 1 & 1 & 1 &1& 1& 11 & 6 \\
\end{array}
\]
As above it follows that $\dim X \leq 6.$ We know that $u_1 \geq
6.$ So, by Lemma \ref{L 3.2} we get the requested constraints in Theorem
\ref{T 2.4}.

Finally we prove Corollary \ref{C 2.3}. Assume that $X$ is not 
arithmetically normal. Then, by Theorem 1.2 of \cite{BS} we know that 
$X$ is a birational projection of a rational normal scroll $Z \subset \mathbb P^{r+1}_K$ 
from a point $p \in \mathbb P^{r+1}_K \setminus Z.$ As $X$ is not a cone, $Z$ is not a cone either and therefore 
$\codim X + 2 = \sum_{i = 1}^k a_i \geq k = \dim X.$ This contradicts the 
assumption of Corollary \ref{C 2.3}. Therefore $X$ is arithmetically 
normal and so $\dim X \leq 6$ by Remark \ref{R 3.3}. 
That is, $\codim X = 3.$ So our claim follows by Theorem \ref{T
2.2}.

For the existence of the samples described in Theorems \ref{T 2.1}
and \ref{T 2.2} we refer to the next section.
\medskip

We close this section with a result on the number of cubics in
a minimal generating set of the defining ideal of a certain
varieties of almost minimal degree.

\begin{lemma} \label{L 3.5} Let $X \subset \mathbb P^r_K$ be a
variety of almost minimal degree with $\codepth A_X = 1$ and $c :=
\codim X \geq 4.$ Then the defining ideal $I_X$ of $X$ is
generated by $\binom{c+1}{2} -2$ quadrics and at most one cubic.
\end{lemma}

\begin{proof} First we reduce the problem to the case in which
$\dim X = 1.$ Let $d = \dim X > 1.$ By an argument of Bertini type
(cf. \cite{J}) we may find generic linear forms
$l_1,\ldots,l_{d-1} \in S_1$ such that $W = X \cap \mathbb
P^{c+1}_K \subset \mathbb P^{c+1}_K :=
\Proj(S/(l_0,\ldots,l_{d-1})S)$ is a non-degenerate integral
variety of almost minimal degree. As $l_0,\ldots, l_{d-1}$ are
chosen generically and $\depth A_X = d$ they form an $A_X$-regular
sequence. Therefore the Betti diagrams of $I_X$ and $I_W$ are the
same. In particular $\codepth A_W = 1.$

So we assume $X \subset \mathbb P^s_K$ with $\dim X = 1$ and $s =
c+1.$ The statement about the number of quadrics is a consequence
of Lemma \ref{L 3.4}. Since $X$ is of almost minimal degree we
know that $I = I_X$ is 3-regular (cf. \ref{L 3.2}). Write $I =
(J,LS)$ with $J = I_2 S$ and with a $K$-vector space $L \subset
S_3$ such that $I_3 = J_3 \oplus L.$ Our aim is to show that
$\dim_K L \leq 1.$

After an appropriate linear coordinate change we may assume that
$x_s \in S_1$ is generic. Let $T := S/x_sS =
K[x_0,\ldots,x_{s-1}].$ Then $R := T/(J,L)T \simeq S/(I,x_sS)$
defines a scheme $Z$ of $s+1$ points in semi uniform position in
$\mathbb P^{s-1}_K.$ The short exact sequence $0 \to A_X(-1)
\stackrel{x_s}{\to} A_X \to R \to 0$ induces an isomorphism
$H^0_{R_+}(R) \simeq K(-2).$ But this means that the vanishing
ideal of $Z$ in $T$ has the form $(J,L,q)T$ with an appropriate
quadric $q \in S_2.$ Since $s \geq 3$ the minimal free resolution
of this ideal has the form
\[
T^{a_2}(-3) \stackrel{\phi}{\to} T^{a_1}(-2) \stackrel{\pi}{\to}
(J,L,q)T \to 0.
\]
This allows us to write $(J,L,q)T = (J,q)T$ and to assume that the
first $a_1-1$ generators of $T^{a_1}(-2)$ are mapped by $\pi$ onto
a $K$-basis of $(JT)_2$ and the last generator is mapped by $\pi$
to $q \cdot 1_T.$ Clearly, $\phi$ is given by a matrix with linear
entries. This shows that $M := JT :_T q \subset T$ is a proper
ideal generated by linear forms.

As $JT \subseteq M$ and as $(J,qT) = IT$ is of height $s-1$ we
must have $s-2 \leq \height M \leq s.$ As $M$ is generated by
linear forms, $(T/M)_1$ is a $K$-vector space of dimension $t  
\in \{0,1,2\}.$ So the graded short exact sequence
\[
0 \to T/M(-2) \to T/JT \to T/(J,q)T \to 0
\]
shows that
\[
\dim_K (IT)_3 = \dim_K ((J,q)T)_3 = \dim_K (JT)_3 + t.
\]
Therefore, we may write $(I,x_s) = (J,L',x_s)$ where $L' \subseteq
L \subset S_3$ is a $K$-vector space of dimension $\leq t.$ As $I$
is a prime ideal and as $x_s \in S_1 \setminus I$ it follows $I =
(J,L'),$ hence $L' = L.$ So, if $\dim_K L' \leq 1,$ we are done.

Otherwise, $\dim_K L' = \dim_K L = 2 = t$ and we may write $I =
(J,k_1,k_2)$ with $k_1, k_2 \in S_3.$ As $\height I = s-1$ it
follows $\height J \geq s-3.$ As $JT \subseteq M$ and as $\height
M = s-2$ we have $\height JT \leq s-2.$ As $x_s$ is a generic
linear form, this means that $\height J \leq s-3$ and hence
$\height J = s-3.$ As $I = (J,k_1,k_2)$ is a prime ideal of height
$s-1 = \height J + 2,$ the ideal $J$ must be prime too.

As $x_s$ is generic and $\height J \leq s-3,$ we may conclude by
Bertini's theorem that $JT \subset T$ defines an integral
subscheme of $\mathbb P^{s-1}_K = \Proj(T).$ So, the saturation
\[
JT :_T \langle T_+ \rangle \subset T \mbox{ of } JT \mbox{ in } T
\]
is a prime ideal of height $s-2.$ As $JT \subseteq M \subset T_+$
and as $M$ is a prime ideal we get $JT :_T \langle T_+ \rangle =
M.$ Therefore
\[
\Proj (T/IT) = \Proj (T/(J,q)T) = \Proj (T/(M,qT))
\]
consists of two points, so that $s + 1 = 2,$ a contradiction. So,
the case $\dim_K L' \geq 2$ does not occur at all.
\end{proof}

The previous Lemma \ref{L 3.5} is inspired by a corresponding
statement for curves of degree $r+2$ in $\mathbb P^r_K$ shown by
the authors (cf. \cite[Lemma (6.4)]{BS1}).

Moreover, if $\codepth A_X \geq 2$ the number of
cubics needed to define $X$ is not bounded by 1 (cf. the examples
in \cite[Section 9]{BS}).

\section{Examples}
In this section we want to confirm the existence of all types of varieties
of almost minimal degree $X \subset P^r_K$ which are described
in the Theorems \ref{T 2.1}, \ref{T 2.2} and \ref{T 2.4}.

First of all we want to show the existence of a Del Pezzo variety
as required by Theorem \ref{T 2.4} (a).

\begin{example} \label{E 4.0} Let $X = \mathbb P^2_K \times
\mathbb P^2_K \subset \mathbb P^8_K$ be the Segre product of two
projective planes. Its defining ideal $I_X$ is generated by the
$2\times 2$-minors of the the following generic $3\times 3$-matrix
\[
\left(%
\begin{array}{ccc}
  x_0 & x_1 & x_2 \\
  x_3 & x_4 & x_5 \\
  x_6 & x_7 & x_8 \\
\end{array}%
\right) \;.
\]
It is easy to see that $\dim X = 4, \codim X = 4$ and $\deg X =
6.$ Therefore, $X$ is a variety of almost minimal degree.
Moreover, $A_X$ is a Cohen-Macaulay and therefore a Gorenstein
ring. An example of dimension 3 is the Segre product $\mathbb P^1_K \times 
\mathbb P^1_K \times \mathbb P^1_K \subset \mathbb P^7_K$ (cf. \cite[(8.11), 6)]{Fu3}). 
Examples of smaller dimensions are obtained by taking generic linear sections. 
\end{example}

In the next examples let us show that the two different Betti
diagrams of statement (b) in \ref{T 2.4} indeed occur. Note that
they require $\codepth A_X = 1$ and $\codim X = 4.$

\begin{example} \label{E 4.7} Consider the rational normal surface
scroll $Z = S(3,3) \subset \mathbb P^7_K.$ Let $P_1 =
(0:0:0:0:1:0:0:1)$ and $P_2 =  (1:0:0:0:0:0:1)$ in $\mathbb
P^7_K.$ Then $P_i \in \mathbb P^7_K \setminus Z, i = 1,2,$ as it
is easily seen. Define $X_i$ to be the projection of $Z$ from
$P_i, i = 1,2.$ Then $\dim X_i = 2$ and $\codepth A_{X_i} = 1$ for
$i= 1,2.$ The Betti diagrams of $I_{X_i}, i = 1,2,$ are those of
Theorem \ref{T 2.4} (b).
\end{example}

Now we construct non-arithmetically Cohen-Macaulay varieties of
almost minimal degree of the type mentioned in Theorem \ref{T
2.2}. To this end we use the possible rational normal scrolls $Z =
S(a_1,\ldots,a_k)$ of the proof of Theorem \ref{T 2.2} which after
appropriate projection furnish the varieties we are looking for.
The construction of the examples corresponding to Theorems \ref{T
2.1} and \ref{T 2.4} follows similarly, and so we skip the details
in these cases.

\begin{example} \label{E 4.2} Let $Z = S(5) \subset \mathbb P^5_K$ denote
the rational normal curve of degree 5. Choose $P \in \mathbb P^5_K
\setminus Z$ a generic point. Then, the
projection $X \subset \mathbb P^4_K$  of $Z$ from $P$ is an example 
of a variety of almost minimal degree with $\dim X = 1$ and
$\codepth A_X = 1.$
\end{example}

Next we want to investigate the case of surfaces.

\begin{example} \label{E 4.3} Let $Z = S(4,1) \subset \mathbb
P^6_K.$ Consider the two points $P_1 = (0:1:0:0:0:0:0)$ and $P_2 =
(0:0:1:0:0:0:0).$ Then $P_i \in \mathbb P^6_K \setminus Z,$ for $i
= 1,2.$ Let $X_i, i = 1,2,$ denote the projection
of $Z$ from $P_i.$ Then $\codepth A_{X_1} = 1$ and $\codepth
A_{X_2} = 2.$ The same type of examples may be produced by
projections of the scroll $S(3,2).$
\end{example}

Our next examples are of dimension 3.

\begin{example} \label{E 4.4} Let $Z = S(3,1,1) \subset \mathbb
P^7_K.$ Consider the points $P_1 = (0:1:0:0:0:0:0:0)$ and $P_2 =
(0:0:0:1:1:0:0:0).$ Then it is easy to see that $P_i \in \mathbb
P^7_K \setminus Z, i = 1,2.$ Let $X_i$ denote the projection of
$Z$ from $P_i.$ Then $\codepth A_{X_1} = 2$ and $\codepth A_{X_1}
= 1.$ The same type of examples may be produced by projections
from the scroll $S(2,2,1).$ 
\end{example}

Now, let us consider the situation of fourfolds.

\begin{example} \label{E 4.5} Consider the scroll $Z = S(2,1,1,1)
\subset \mathbb P^8_K.$ Let $P_1 = (0:1:0:0:0:0:0:0:0)$ and $P_2 =
(0:0:0:0:0:0:1:1:0).$ Then $P_i \in \mathbb P^8_K \setminus Z, i =
1,2.$ Let $X_i \subset \mathbb P^7$ denote the projection of $Z$
from $P_i, i = 1,2.$ Then $\codepth A_{X_1} = 2,$ while $\codepth
A_{X_2} = 1.$
\end{example}

Finally let us consider the case where $\dim X = 5.$

\begin{example} \label{E 4.6} Let $Z = S(1,1,1,1,1) \subset
\mathbb P^9_K$ be the Segre variety. Then $P \in \mathbb P^9_K
\setminus Z$ for the point $P = (0:1:1:0:0:0:0:0:0:0).$ Let $X
\subset \mathbb P^8_K$ denote the projection of $Z$ from $P.$ Then
$\depth A_X = 4,$ and therefore $\codepth A_X = 2.$ Finally 
observe that $\codepth A_X = 1$ is impossible if $\dim X =
5,$ as $\depth A_X  \leq 4$ (cf. Lemma \ref{L 3.1}).
\end{example}

The Examples \ref{E 4.2} -- \ref{E 4.6} provide the existence of
the samples claimed by Theorem \ref{T 2.2}. Similar constructions
provide varieties as mentioned in Theorem \ref{T 2.1} and \ref{T
2.4}.

In the final examples, we will show that in higher codimension, the
shape of the Betti diagram of $I_X$ for a variety $X$ of minimal 
degree may vary in a much stronger way: In fact the "beginning of
the Betti diagrams" may be rather different from each other.

\begin{example} \label{E 4.1} Let $Z = S(8) \subset \mathbb P^8_K$ denote
the rational normal curve of degree 8. Let $P_1 =
(0:0:0:0:0:0:1:0:0),P_2 = (0:0:0:0:0:1:0:0:0)$ and $P_3 =
(0:0:0:0:1:0:0:0:0).$ Then $P_i \in \mathbb P^8_K \setminus Z$ for
$i = 1,2,3.$ Let $X_i \subset \mathbb P^7_K, i = 1,2,3,$ denote
the projection of $Z$ from $P_i.$ Then the Betti diagrams of
$I_{X_i}, i = 1,2,3,$ resp. have the form:
\[
\begin{tabular}{|c||c|ccccccc|}
  \hline
   $i$ & & 1 & 2 & 3 & 4 & 5 & 6 & 7\\
  \hline
 1 & 1 & 19 & 58 & 75 &  44 & 5 & 0 & 0 \\
  & 2 & 1  & 6  & 15 &  20 & 21& 8 & 1 \\
  \hline
 2 & 1 & 19 & 57 & 70 &  34 & 5 & 0 & 0 \\
  & 2 & 0  & 1 & 5 &  20 & 21& 8 & 1 \\
  \hline
 3 & 1 & 19 & 57 & 69 &  34 & 5 & 0 & 0 \\
  & 2 & 0  & 0 & 5 &  20 & 21& 8 & 1 \\
  \hline
\end{tabular}
\]
In all three cases $\codim X_i = 4$ and $\codepth A_{X_i} =1.$ 
Remember that the number of cubics in the defining ideals
is bounded by 1 (cf. Lemma \ref{L 3.5}). It follows that $X_1$ is
contained in the scroll $S(5,1),$ while $X_2$ is contained in the
scroll $S(4,2)$ and $X_3$ is contained in the scroll $S(3,3).$
\end{example}

In view of the Example \ref{E 4.1} and corresponding examples in
higher dimensions one might expect that the type of the rational
normal scroll $Y,$ that contains the variety $X$ of almost minimal degree 
as a one codimensional subvariety,  
determines the Betti diagram "near the beginning of the
resolution". In small codimensions, the different types of these 
scrolls $Y$ are much more limited than imposed by Theorems \ref{T
2.1}, \ref{T 2.2} and \ref{T 2.4}. It seems rather challenging to
understand the r\^ole of the scrolls $Y$ for the beginning of the
minimal free resolution of $I_X.$

Moreover the examples in \ref{E 4.1} show that the estimates for
the Betti numbers given in Lemma \ref{L 3.2} near the beginning of 
the Betti diagram are fairly weak.

\begin{remark} \label{R 4.7}
To compute the Betti diagrams and hence the arithmetic depths of
the above examples, we have made use of the computer algebra
system {\sc Singular} (cf. \cite{G}). Moreover, there is in
preparation a conceptual approach for the computation of the
$\depth A_X$ in terms of the center of the projection and the
secant variety of $S(a_1,\ldots,a_k) \subset \mathbb P^{r+1}_K.$
\end{remark}


\begin{thebibliography}{9999}

\bibitem[1]{BS1} {\sc Brodmann, M., Schenzel, P.}: {\em Curves of
   degree $r+2$ in $\mathbb P^r$: Cohomological, geometric and
   homological aspects,} J. Algebra 242 (2001), 577 - 623.


\bibitem[2]{BS} {\sc Brodmann, M., Schenzel, P.}: {\em Arithmetic
   properties of projective varieties of almost minimal degree,}
   to appear in J. Algebraic Geom.

\bibitem[3]{Fu1} {\sc Fujita, T.}: {\em Classification of projective
   varieties of $\Delta $-genus one}. Proceedings of the Japan
   Academy of Science, Ser. A Math. Sci 58 (1982), 113 - 116.

\bibitem[4]{Fu2} {\sc Fujita, T.}: {\em Projective varieties of
   $\Delta $-genus one}. In Algebraic and Topological Theories --
   To the Memory of Dr. Takehiko Miyata (Kinokuniya, Tokyo, 1985),
   149 - 175.

\bibitem[5]{Fu3} {\sc Fujita, T.}: {\em Classification theories of
   polarized varieties}, London Mathematical Society Lecture
   Notes Series 155, Cambridge University Press, 1990.

\bibitem[6]{E} {\sc Eisenbud, D.}: {\em The Geometry of Syzygies
   -- A second course in Commutative Algebra and Algebraic Geometry,}
   Graduate Texts in Mathematics, Vol. 229, Springer-Verlag, New York,
   2005.

\bibitem[7]{EH} {\sc Eisenbud, D., Harris, J.}: {\em On varieties of
   minimal degree (a centennial account)}. In Proceedings of Symposion
   of Pure Mathematics Vol 46, pp 3 - 13, American Mathematical Society,
   Providence, 1987.

\bibitem[8]{G} {\sc Greuel, G.M., Pfister, G. et al}: {\em Singular
   $3.0$, A Computer Algebra System for Polynomial Computations}. Center for Computer
   Algebra, University of Kaiserslautern (2005)
   (http://www.singular.uni-kl.de).

\bibitem[9]{H} {\sc Harris, J.}: {\em Algebraic Geometry: A First Course}.
   Graduate Texts in Mathematics, Vol. 133, Springer-Verlag, New York,
   1992.

\bibitem[10]{HSV} {\sc Hoa, T., St\"uckrad, J., Vogel, W.}: {\em Towards
   a structure theory of projective varieties of degree $=$ codimension
   $+ 2$}. J. of Pure and Appl. Algebra 71 (1991), 203 - 231.

\bibitem[11]{J} {\sc Jouanolou, J.P.}: {\em Th\'eor\`emes de Bertini et
   applications}. Progress in Mathematics, Vol. 42, Birkh\"auser,
   Basel, 1983.

\bibitem[12]{Na} {\sc Nagel,U}: {\em Minimal free resolutions of
   projective subschemes of small degree,} Preprint, 2005, 
   (http://arxiv.org/math.AG/0505347).

\bibitem[13]{S} {\sc Schenzel, P.}: {\em \"Uber die freien Aufl\"osungen
   extremaler Cohen-Macaulay-Ringe,} J. Algebra 64, 93-101 (1980).

\end{thebibliography}
\end{document}